\documentclass[conference, 10pt, twocolumn]{ieeeconf}       

\def\BibTeX{{\rm B\kern-.05em{\sc i\kern-.025em b}\kern-.08em
    T\kern-.1667em\lower.7ex\hbox{E}\kern-.125emX}}
    
\usepackage[left=54pt, right=54pt,bottom=54pt, top=54pt]{geometry}

\usepackage{amsmath,mathrsfs,amsfonts,amssymb,graphicx,epsfig,mathtools}
 \usepackage{amsthm}
\usepackage{subcaption}
\usepackage{color,multirow,rotating}
\usepackage{algorithm,algpseudocode,algorithmicx}
\usepackage{cite,url,framed,bm,balance}

\usepackage{dsfont}

\usepackage{varwidth}
\usepackage{cleveref}

\usepackage{tikz}
\usetikzlibrary{calc,trees,positioning,arrows,chains,shapes.geometric,decorations.pathreplacing,decorations.pathmorphing,shapes, matrix,shapes.symbols}

\DeclareMathOperator*{\arginf}{arg\:inf}

\newcommand{\differential}{{\rm{d}}}

\newcommand{\bsin}{\bm{{\sin}}}
\newcommand{\bcos}{\bm{{\cos}}}
\newcommand{\prox}{{\mathrm{prox}}}
\newcommand{\RNum}[1]{\uppercase\expandafter{\romannumeral #1\relax}}

\allowdisplaybreaks

\title{\LARGE\textbf{
Schr\"{o}dinger Meets Kuramoto via Feynman-Kac: Minimum Effort Distribution Steering for Noisy Nonuniform Kuramoto Oscillators}
}

\author{Iman Nodozi and Abhishek Halder
\thanks{Iman Nodozi and Abhishek Halder are with the University of California, Santa Cruz, CA 95064, USA,
        {\tt\small{\{inodozi,ahalder\}@ucsc.edu}}. This work is partially supported by NSF grants 1923278, 2112755.%
}}

\IEEEoverridecommandlockouts
\begin{document}

\maketitle
\pagenumbering{gobble}

\begin{abstract}
We formulate and solve the problem of finite horizon minimum control effort steering of the state probability distribution between prescribed endpoint joints for a finite population of networked noisy nonuniform Kuramoto oscillators. We consider both the first and second order stochastic Kuramoto models. For numerical solution of the associated stochastic optimal control, we propose combining certain measure-valued proximal recursions and the Feynman-Kac path integral computation. We illustrate the proposed framework via numerical examples.
\end{abstract}


\section{Introduction}\label{sec:introduction}
We consider the controlled sample path dynamics for a population of $n$ \emph{first order} Kuramoto oscillators, given by the It\^{o} stochastic differential equations (SDEs)
$$\differential\theta_{i} = \left(-\frac{\partial V}{\partial\theta_i} + v_{i}\right)\differential t + \sqrt{2}\sigma_{i}\:\differential w_{i}, \; i\in[n]:=\{1,2,\hdots,n\},$$
where $V(\theta_1,\hdots,\theta_n)$ is a given smooth potential, the angular variable $\theta_{i}\in [0,2\pi)$ is the state, $v_{i}$ is the control input, $\sigma_i>0$ is the noise strength, and $w_{i}$ is the standard (scalar) Wiener process noise for the $i$th oscillator. Defining re-scaled input $u_{i}:=v_{i}/\sigma_{i}$, we write this dynamics in vector form:
\begin{align}
\differential\bm{\theta} = \left(-\nabla_{\bm{\theta}} V(\bm{\theta}) + \bm{S} \bm{u}\right)\differential t + \sqrt{2} \bm{S}\:\differential\bm{w}
\label{FirstOrderKuramotoSDE}    
\end{align}
where $\bm{S}:={\rm{diag}}\left(\sigma_{1},\ldots,\sigma_{n}\right)\succ\bm{0}$, $\bm{\theta}:=(\theta_1,\hdots,\theta_n)^{\top}$, $\bm{u}:=(u_1,\hdots,u_n)^{\top}$, and $\bm{w}:=(w_1,\hdots,w_n)^{\top}$ is the standard Wiener process in $n$ dimensions. For the first order Kuramoto model \eqref{FirstOrderKuramotoSDE}, the state space is the $n$-torus $\mathbb{T}^{n} \equiv [0,2\pi)^{n}$, and the potential
\begin{align}
    V(\bm{\theta}):=\!\!\sum\limits_{\stackrel{i<j}{i,j\in[n]}}\!\! k_{ij} (1-\cos(\theta_i-\theta_j-\varphi_{ij}))-\sum\limits_{i=1}^{n}P_i \theta_i,
    \label{V_1}
\end{align}  
wherein the parameters $P_{i}>0$. For $i\neq j$, the coupling coefficients $k_{ij}=k_{ji}\geq 0$ (and not all $k_{ij}=0$), $k_{ii}\equiv 0$. Likewise, for $i\neq j$, the phase shift $\varphi_{ij}=\varphi_{ji}\in[0,\frac{\pi}{2})$, and $\varphi_{ii}\equiv 0$. 

We also consider the controlled sample path dynamics for a population of $n$ \emph{second order} Kuramoto oscillators, given by the second order Langevin equations
\begin{align*}
m_{i}\ddot{\theta}_{i} + \gamma_{i}\dot{\theta}_{i} = &-\frac{\partial V}{\partial\theta_i} + v_i +  \sqrt{2}\sigma_{i}\times\text{SGWN},    
\end{align*}
where $i\in[n]$, SGWN denotes standard Gaussian white noise, and $V(\cdot)$ is given by \eqref{V_1}. Letting $u_{i}:=v_{i}/\sigma_{i}$ as before, we rewrite this second order dynamics as the vector It\^{o} SDE
\begin{align}
\begin{pmatrix}
\differential\bm{\theta}\\
\differential\bm{\omega}
\end{pmatrix}\! =& \!\begin{pmatrix}
\bm{\omega}\\
- \bm{M}^{-1}\nabla_{\bm{\theta}} V(\bm{\theta}) - \bm{M}^{-1} \bm{\Gamma}\bm{\omega} + \bm{M}^{-1} \bm{S} \bm{u}
\end{pmatrix}\!\differential t~ \!\nonumber\\
&  + \begin{pmatrix}
\bm{0}_{n\times 1}\\
\sqrt{2}\bm{M}^{-1}\bm{S}\:\differential\bm{w}
\end{pmatrix}
\label{SecondOrderKuramotoSDE}    
\end{align}
where $\bm{\omega}:=(\dot{\theta}_{1},\hdots,\dot{\theta}_{n})^{\top}$, $\bm{M}:= {\rm{diag}}\left(m_{1},\ldots,m_{n}\right)\succ\bm{0}$,
$\bm{\Gamma} := {\rm{diag}}\left(\gamma_{1},\ldots,\gamma_{n}\right)\succ\bm{0}$, and $\bm{0}_{n\times 1}$ denotes the $n\times 1$ vector of zeros. For the second order Kuramoto model, the state space is the product of cylinders $\mathbb{T}^{n}\times\mathbb{R}^{n}$.

In this paper, we address the following problem:\\
\noindent synthesize minimum effort control policy $\bm{u}$ that transfers the stochastic state of \eqref{FirstOrderKuramotoSDE} or \eqref{SecondOrderKuramotoSDE} from a prescribed initial to a prescribed terminal joint probability distribution over a given finite time horizon, say $t\in[0,T]$. 
 
This fits in the research theme of designing state feedback for dynamically reshaping (as opposed to simply mitigating) uncertainties \cite{brockett2012notes,chen2021controlling} subject to networked Kuramoto oscillator dynamics. As such, both first and second order Kuramoto oscillator models are ubiquitous across physical, biological and engineering systems, see e.g., \cite{acebron2005kuramoto,dorfler2014synchronization}.

Notice that while the uncontrolled dynamics in \eqref{FirstOrderKuramotoSDE} has gradient drift, the same in \eqref{SecondOrderKuramotoSDE} has mixed conservative-dissipative drift. A consequence is that unlike \eqref{FirstOrderKuramotoSDE}, the stochastic process induced by \eqref{SecondOrderKuramotoSDE}, is \emph{not} reversible and its infinitesimal generator is hypoelliptic \cite{hormander1967hypoelliptic}. This makes the analysis and feedback synthesis for \eqref{SecondOrderKuramotoSDE} even more challenging than \eqref{FirstOrderKuramotoSDE}.

\subsubsection*{Related literature and novelty of this work} While there exists a significant literature on the dynamics  and control of Kuramoto oscillators in general \cite{strogatz2000kuramoto,jadbabaie2004stability,chopra2009exponential,yin2011synchronization,dorfler2012synchronization,sahyoun2015optimal,li2021unified}, the \emph{stochastic control} of Kuramoto oscillators remains under-investigated. Ref. \cite{wu2021global} considered global asymptotic phase agreement and frequency synchronization in almost sure sense.

In the physics literature, several studies \cite{martens2009exact,benedetto2015complete,bertini2010dynamical} analyze the distributional dynamics associated with the Kuramoto oscillators. However, these studies consider the \emph{univariate} distributional dynamics arising from the \emph{mean-field limit}, i.e., by abstracting the dynamical interaction in the \emph{infinite population} ($n\rightarrow \infty$) regime. In comparison, the perspective and approach taken in this paper are significantly different because we focus on the dynamics of \emph{joint} probability distribution supported over the states of a \emph{finite population} of oscillators. This is particularly relevant for engineering applications such as power systems, where a network of finitely many generators (often modeled as second order nonuniform Kuramoto oscillators) and loads (often modeled as first order nonuniform Kuramoto oscillators) interact together with their controlled stochastic dynamics, see e.g., \cite{dorfler2012synchronization}. Well-known techniques such as the Kron reduction \cite{dorfler2012kron} allow transcribing such networked system in the form \eqref{SecondOrderKuramotoSDE} with all-to-all connection topology. Despite the engineering relevance, research on the multivariate distributional dynamics for a finite population of nonuniform Kuramoto oscillators is scant. 

From a methodological standpoint, we recast the problem of minimum effort feedback steering of distributions subject to \eqref{FirstOrderKuramotoSDE} or \eqref{SecondOrderKuramotoSDE}, as an instance of generalized Schr\"{o}dinger bridge problem -- a topic undergoing rapid development \cite{chen2021stochastic,chen2016relation} in the systems-control community. In \cite{chen2015fast}, a similar approach was taken to realize feedback steering toward the invariant distribution of an uncontrolled oscillator dynamics. Building on our prior work \cite{caluya2021wasserstein}, here we focus on finite horizon steering between two arbitrary compactly supported joint state probability distributions subject to \eqref{FirstOrderKuramotoSDE} or \eqref{SecondOrderKuramotoSDE}. However, for our controlled Kuramoto dynamics, it will turn out that the algorithmic approach proposed in \cite{caluya2021wasserstein} will no longer apply and we will introduce new ideas for the same.


\subsubsection*{Notations} We use boldfaced capital letters for matrices, and boldfaced small letters for vectors. The symbol $\mathbb{E}_{\mu^{\bm{u}}}\left[\cdot\right]$ denotes the mathematical expectation w.r.t. the controlled joint state probability measure $\mu^{\bm{u}}$, that is, $\mathbb{E}_{\mu^{\bm{u}}}\left[\cdot\right] := \int (\cdot)\:\differential\mu^{\bm{u}}$. The superscript $\bm{u}$ in $\mu^{\bm{u}}$ indicates that the joint measure depends on the choice of control $\bm{u}$. For the controlled dynamics \eqref{FirstOrderKuramotoSDE}, the measure $\mu^{\bm{u}}$ is supported over the state space $\mathbb{T}^{n}$. Likewise, for \eqref{SecondOrderKuramotoSDE}, the measure $\mu^{\bm{u}}$ is supported over $\mathbb{T}^{n}\times\mathbb{R}^{n}$. The symbol $\sim$ is used as a shorthand for ``follows the probability distribution". The notations $\nabla$, $\Delta$, $\textbf{Hess}(\cdot)$, $\langle\cdot,\cdot\rangle$, $\otimes$, $\bm{I}_{n}$ respectively denote the Euclidean gradient, Laplacian, Hessian, the Euclidean inner product (Frobenius inner product for matricial arguments), the Kronecker product, and the $n\times n$ identity matrix.

\subsubsection*{Organization} The outline of this paper is as follows. Sec. \ref{sec:problem} details the problem formulation. The existence and uniqueness of its solution are discussed in Sec. \ref{sec:existenceuniqueness}. In Sec. \ref{sec:properties}, we detail how the optimal solutions can be recovered from the so-called Schr\"{o}dinger factors which in turn, solve a nonlinearly boundary-coupled system of linear PDEs. We derive these systems for both the first and second order Kuramoto oscillators. Sec. \ref{sec:algo} summarizes the proposed combination of proximal and Feynman-Kac algorithms for solving the respective boundary-coupled systems, followed by numerical simulations in Sec. \ref{sec:numericalsimulation}. Concluding remarks are provided in Sec. \ref{sec:conclusions}.


\section{The Optimal Distribution Steering Problem}\label{sec:problem}
\subsubsection{Formulation}
We consider a stochastic optimal control problem over prescribed time horizon $[0,T]$, given by
\begin{align}
\underset{\bm{u}\in\mathcal{U}}{\inf}\:\mathbb{E}_{\mu^{\bm{u}}}\left[\int_{0}^{T}\|\bm{u}\|_{2}^{2}\:\differential t\right]
\label{SBPobjective}    
\end{align}
subject to \emph{either}
\begin{align*}
\eqref{FirstOrderKuramotoSDE},\; \bm{\theta}(t=0)\sim \mu_{0}\;\text{(given)}, \; \bm{\theta}(t=T)\sim \mu_{T}\;\text{(given)},
\end{align*}
\emph{or} 
\begin{align*}
\eqref{SecondOrderKuramotoSDE},\; \begin{pmatrix}\bm{\theta}(t=0)\\
\bm{\omega}(t=0)\end{pmatrix}\sim \mu_{0}\;\text{(given)}, \; \begin{pmatrix}\bm{\theta}(t=T)\\
\bm{\omega}(t=T)\end{pmatrix}\sim \mu_{T}\;\text{(given)},
\end{align*}
where $\mu_{0},\mu_{T}$ denote the joint state probability measures at $t=0$ and $t=T$, respectively. In \eqref{SBPobjective}, the feasible set $\mathcal{U}$ comprises of the finite energy Markovian state and time dependent input policies over the time horizon $[0,T]$.

Assuming the absolute continuity of the joint probability measure $\mu^{\bm{u}}$ for all times, we write $\differential \mu^{\bm{u}} (\bm{x},t)=\rho^{\bm{u}}(\bm{x},t) \differential\bm{x}$ and hereafter consider the associated joint PDF $\rho^{\bm{u}}(\bm{x},t)$. Problem \eqref{SBPobjective} can then be recast as
\begin{align}
\underset{(\rho^{\bm{u}},\bm{u})}{\inf}\:\int_{0}^{T}\int_{\mathcal{X}}\|\bm{u}(\bm{x},t)\|_{2}^{2}\:\rho^{\bm{u}}(\bm{x},t)\:\differential \bm{x} \:\differential t
\label{FPKobjective}    
\end{align}
subject to either
\begin{subequations}
	\begin{align}
		&\frac{\partial \rho^{\bm{u}}}{\partial t} =-\nabla_{\bm{\theta}} \cdot(\rho^{\bm{u}} ( \bm{S} \bm{u}-\nabla_{\bm{\theta}} V))+\langle \bm{D}, \textbf{Hess}(\rho^{\bm{u}})\rangle, \label{FPPDEFirstOrder} \\
	 & \begin{aligned}
	   	\!\!\!\!\text{or}\;&\frac{\partial \rho^{\bm{u}}}{\partial t} =\nabla_{\bm{\omega}} \cdot(\rho^{\bm{u}} \left( \bm{M}^{-1}\nabla_{\bm{\theta}} V(\bm{\theta}) + \bm{M}^{-1} \bm{\Gamma}\bm{\omega} - \bm{M}^{-1} \bm{S} \bm{u}\right. \\
	   	&\left.+\bm{M}^{-1}\bm{D}\bm{M}^{-1} \nabla_{\bm{\omega}}\log\rho^{\bm{u}} \right) -\langle \bm{\omega},\nabla_{\bm{\theta}} \rho^{\bm{u}} \rangle,\!\!\! \label{FPPDESecondOrder}
	   	\end{aligned}
	\end{align}
	\label{FPK_PDEs}
\end{subequations}
where the diffusion matrix $\bm{D}:=\bm{S} \bm{S}^{\top}$, and $\rho^{\bm{u}}(\bm{x},t=0)=\rho_0$ (given), $\rho^{\bm{u}}(\bm{x},t=T)=\rho_T$ (given).

For the first order Kuramoto oscillators, we have $\bm{x}:=\bm{\theta}$, $\mathcal{X}:=\mathbb{T}^n$, and for the second order Kuramoto oscillators, we have $\bm{x}:=(\bm{\theta}, \bm{\omega})^\top$, $ \mathcal{X}:=\mathbb{T}^n \times \mathbb{R}^n$. The constraints \eqref{FPPDEFirstOrder} and \eqref{FPPDESecondOrder} are the controlled Fokker-Planck-Kolmogorov (FPK) forward PDEs corresponding to \eqref{FirstOrderKuramotoSDE} and \eqref{SecondOrderKuramotoSDE}, respectively.

\subsubsection{Endpoint PDFs}
In this work, we suppose that the endpoint joint PDFs $\rho_{0}, \rho_{T}$ are supported on compact subsets of $\mathcal{X}$. For instance, when $\mathcal{X}=\mathbb{T}^{n}$, one may model $\rho_{0},\rho_{T}$ as multivariate von Mises PDFs \cite{mardia2008multivariate,mardia2014some} supported on $\mathbb{T}^{n}$:
\begin{align}
    &\!\rho_{k}(\bm{\theta}) =\frac{1}{Z_{k}} \exp\bigl(\langle\bm{\kappa}_{k},\bcos(\bm{\theta}-\bm{m}_{k})\rangle \nonumber\\
    & +\frac{1}{2}\langle\bsin(\bm{\theta}-\bm{m}_{k}),\bm{\Lambda}_{k}\bsin(\bm{\theta}-\bm{m}_{k})\rangle\!\bigr), \quad k\in\{0,T\},
    \label{endpointPDFs}
\end{align}
where the parameters are the mean vectors $\bm{m}_{0},\bm{m}_{T}\in \mathbb{T}^{n}$, the concentration vectors $\bm{\kappa}_{0},\bm{\kappa}_{T}\in \mathbb{R}_{\geq 0}^{n}$, and  $\bm{\Lambda}_{0},\bm{\Lambda}_{T}\in \mathbb{S}^{n}$ ($n\times n$ real symmetric matrices) having zero diagonal entries. In \eqref{endpointPDFs}, $\bm{{\sin}}(\cdot)$ and $\bm{{\cos}}(\cdot)$ denote the elementwise sines and cosines, respectively. The normalization constants $Z_{0},Z_{T}$ in \eqref{endpointPDFs} depend on the respective concentration vector and symmetric matrix parameters.

The nonnegative entries of the concentration vectors $\bm{\kappa}_{0},\bm{\kappa}_{T}$ admit a natural interpretation: zero concentration vectors represent uniform distribution over $\mathbb{T}^{n}$. Large positive entries promote a higher concentration around the corresponding mean components. When $\bm{\Lambda}$ is a zero matrix, then multivariate von Mises PDF can be written as the product of univariate von Mises PDFs, see e.g., \cite[Ch. 3]{mardia2009directional}.

When $\mathcal{X}=\mathbb{T}^{n}\times\mathbb{R}^{n}$, we suppose that for $k\in\{0,T\}$, the $\bm{\omega}$ marginals of $\rho_{k}$ have compact supports $\Omega_{k}\subset\mathbb{R}^{n}$, and thus the joints $\rho_{k}$ are supported on compact subsets of $\mathcal{X}$.


\section{Existence and Uniqueness of Solution}\label{sec:existenceuniqueness}
\subsubsection{First order case}\label{subsec:FirstOrderExistenceUniqeness}
From \eqref{V_1}, we observe that $V\in C^{2}(\mathbb{T}^n)$, which allows us \cite[Ch. 1.2]{stroock2008partial} to conclude that the transition probability kernels associated with \eqref{FirstOrderKuramotoSDE} remain continuous for all $t \geq 0$. Furthermore, the endpoint PDFs having compact supports imply $\rho_{0}, \rho_{T}$ are positive over their respective supports. Thus, following \cite[Appendix E]{caluya2021wasserstein}, the transition probability kernels associated with \eqref{FirstOrderKuramotoSDE} also remain positive for all $t \in [0,T]$.

The continuity and positivity of the transition probability kernels associated with \eqref{FirstOrderKuramotoSDE}, together guarantee \cite[Sec. 10]{wakolbinger1990schrodinger}, \cite[Thm. 3.2]{jamison1974reciprocal} the existence-uniqueness for the solution of the variational problem \eqref{FPKobjective} subject to \eqref{FPPDEFirstOrder} and the endpoint PDF constraints.

\subsubsection{Second order case}\label{subsec:SecondOrderExistenceUniqeness}
That the transition probability kernels remain positive, is ensured per the compactness assumption of the endpoint joint PDFs' supports together with the maximum principle for parabolic PDEs. 

Showing that the transition probability kernels also remain continuous for all times, in this case, reduces to showing three conditions: (i) $V\in C^{2}\left(\mathbb{T}^{n}\right)$, (ii) $\inf\:V > -\infty$, and (iii) uniform boundedness of the Hessian: $\|\textbf{Hess}(V)\|_{2} \leq c$ for some $c>0$ that does not depend on $\bm{\theta}$; see e.g., \cite[Theorem 7]{villani2009hypocoercivity}, \cite[Theorem 5]{markou2014fokker}. The satisfaction of the conditions (i)-(ii) are immediate. For condition (iii), notice that the induced 2-norm of $\textbf{Hess}(V)$ is upper bounded by
$$\sum_{i<j}k_{ij}\cos\left(\theta_{i}-\theta_{j}-\varphi_{ij}\right) \leq \sum_{i<j}k_{ij}.$$
Since $k_{ij}\geq 0$ for all $i,j\in[n]$, and there exists $i,j\in[n]$ such that $k_{ij}>0$, therefore, (iii) also holds.

As in the first order case, the continuity and positivity of the transition probability kernels, together guarantee the existence-uniqueness of the solution of \eqref{FPKobjective} subject to \eqref{FPPDESecondOrder} and the endpoint PDF constraints.

In the following Section, we express the solutions of \eqref{FPKobjective} in terms of the so-called Schr\"{o}dinger factors for both first and second order controlled Kuramoto dynamics.


\section{Optimal Solutions and Schr\"{o}dinger Factors}\label{sec:properties}
\subsection{First Order Case}
Since $\bm{S}$ is not identity, the strengths of the process noise acting along the components of \eqref{FirstOrderKuramotoSDE} are nonuniform. To account this anisotropic noise, we consider an invertible linear map $\bm{\theta}\mapsto \bm{\xi}:=\bm{S}^{-1} \bm{\theta}$, which by It$\hat{\text{o}}$'s Lemma \cite[Ch.4.2]{oksendal2013stochastic}, results in the following SDE for the transformed state vector $\bm{\xi}$:
\begin{equation}
    \begin{aligned}
     \differential\bm{\xi} = \left(\bm{u}- \bm{\Upsilon} \nabla_{\bm{\xi}} \Tilde{V}(\bm{\xi}) \right)\differential t + \sqrt{2}\:\differential\bm{w}
    \end{aligned}
    \label{SDE_xi}
\end{equation}
where the matrix $\bm{\Upsilon}\!:=\!\left(\prod_{i=1}^{n}\sigma_{i}^{2}\right)\!\bm{S}^{-2}\!=\!{\rm{diag}}\!\left(\!\prod_{j\neq i} \sigma_{j}^{2}\!\right)\!\succ\!\bm{0}$, and the potential
{\small{\begin{equation*}
         \begin{aligned}
 \Tilde{V}(\bm{\xi}):=&\left( \sum\limits_{i<j}k_{ij}\left(1-\cos(\sigma_{i}\xi_{i}-\sigma_{j}\xi_{j}-\varphi_{ij})\right)-\right.\\
&\left.\sum\limits_{i=1}^{n} \sigma_{i}P_{i}\xi_{i} \right) \bigg/ \left(\prod\limits_{i=1}^{n}\sigma_{i}^{2}\right).
         \end{aligned}
         \end{equation*}}}

In this new state coordinate, the problem \eqref{FPKobjective} subject to \eqref{FPPDEFirstOrder} and the endpoint PDF constraints, takes the form
\begin{subequations}
\begin{align}
&\underset{(\tilde{\rho}^{\bm{u}},\bm{u})}{\inf}\:\int_{0}^{T}\int_{\mathcal{X}}\|\bm{u}(\bm{\xi},t)\|_{2}^{2}\:\tilde{\rho}^{\bm{u}}(\bm{\xi},t)\:\differential \bm{\xi}\:\differential t\\
&	\frac{\partial \tilde{\rho}^{\bm{u}}}{\partial t} =-\nabla_{\bm{\xi}} \cdot(\tilde{\rho}^{\bm{u}} (\bm{u}-\bm{\Upsilon} \nabla_{\bm{\xi}} \tilde{V}))+\Delta_{\bm{\xi}} \tilde{\rho}^{\bm{u}}, \\
&\begin{aligned}
&\tilde{\rho}^{\bm{u}}(\bm{\xi},0)=\rho_{0}(\bm{S}\bm{\xi})\!\!\left(\prod_{i=1}^{n} \sigma_i\right),  \; \tilde{\rho}^{\bm{u}}(\bm{\xi},T)=\rho_{T}(\bm{S}\bm{\xi})\!\!\left(\prod_{i=1}^{n} \sigma_i\right).
\label{endpointPDFsNewCordinate}
\end{aligned}
\end{align}
\label{opt_first_order}
\end{subequations}
Applying Proposition 1 and Theorem 2 of \cite{caluya2021wasserstein} to \eqref{opt_first_order}, we derive a boundary-coupled system of linear PDEs for the function pair $(\varphi(t,\bm{\xi}),\hat{\varphi}(t,\bm{\xi}))$, given by
\begin{subequations}
\begin{align}
&\frac{\partial \hat{\varphi}}{\partial t}= \nabla_{\bm{\xi}}\cdot (\hat{\varphi} \bm{\Upsilon} \nabla_{\bm{\xi}} \Tilde{V})+\Delta_{\bm{\xi}} \hat{\varphi},~~~\label{Forward_PDE_First_Order}\\
&\frac{\partial \varphi}{\partial t}=\langle \nabla_{\bm{\xi}}\varphi, \bm{\Upsilon} \nabla_{\bm{\xi}} \Tilde{V}\rangle-\Delta_{\bm{\xi}} \varphi,\label{Backward_PDE_First_Order}\\
&\hat{\varphi}_{0}(\bm{\xi})\varphi_{0}(\bm{\xi}) = \tilde{\rho}^{\bm{u}}(\bm{\xi},0) = \rho_{0}(\bm{S}\bm{\xi})\!\!\left(\prod_{i=1}^{n} \sigma_i\right),\label{CoupledBoundaryCondFirstOrder1}\\
&\hat{\varphi}_{T}(\bm{\xi})\varphi_{T}(\bm{\xi}) = \tilde{\rho}^{\bm{u}}(\bm{\xi},T)=\rho_{T}(\bm{S}\bm{\xi})\!\!\left(\prod_{i=1}^{n} \sigma_i\right), \label{CoupledBoundaryCondFirstOrder2}
\end{align}
\label{Coupled_PDE_first_order}
\end{subequations}
whose solution recovers the optimal decision variables $(\tilde{\rho}^{\text{opt}},\bm{u}^{\text{opt}})$ for problem \eqref{opt_first_order} via the mapping
\begin{align}
 \tilde{\rho}^{\text{opt}}(\bm{\xi},t)\!=\! \hat{\varphi}(\bm{\xi},t)\varphi(\bm{\xi},t), \; \bm{u}^{\text{opt}}(\bm{\xi},t)\!=\!\nabla_{\bm{\xi}} \log \varphi (\bm{\xi},t).
  \label{opt_pair}
\end{align}
We refer to the function pair $(\varphi,\hat{\varphi})$ as the \emph{Schr\"{o}dinger factors}, so named since their product gives $\tilde{\rho}^{\text{opt}}$ at all times, i.e., $(\varphi,\hat{\varphi})$ comprise a factorization of $\tilde{\rho}^{\text{opt}}$. The optimally controlled joint state PDF $\rho^{\text{opt}}$ for \eqref{FPKobjective} is then obtained as $\rho^{\text{opt}}(\bm{\theta},t) = \tilde{\rho}^{\text{opt}}\left(\bm{S}^{-1}\bm{\theta},t\right)/\left(\prod_{i=1}^{n}\sigma_{i}\right)$. The optimal control in original coordinates is $\bm{S}\nabla_{\bm{\theta}}\log\varphi(\bm{S}^{-1}\bm{\theta},t)$.

Now the matter boils down to solving \eqref{Coupled_PDE_first_order}. For notational ease, let $\hat{\varphi}_{0}:=\hat{\varphi}(\bm{\xi},0)$, $\hat{\varphi}_{T}:=\hat{\varphi}(\bm{\xi},T)$, $\varphi_{0}:=\varphi(\bm{\xi},0)$, and $\varphi_{T}:=\varphi(\bm{\xi},T)$. Notice that \eqref{Forward_PDE_First_Order}-\eqref{Backward_PDE_First_Order} are the \emph{uncontrolled} forward and backward Kolmogorov PDEs, respectively, associated with \eqref{SDE_xi}. Since \eqref{Forward_PDE_First_Order}-\eqref{Backward_PDE_First_Order} are equation-level-decoupled, the system \eqref{Coupled_PDE_first_order} can be seen as a nonlinear fixed point map for the pair $\left(\hat{\varphi}_0,\varphi_{T}\right)$ that is known \cite{chen2016entropic} to be contractive w.r.t. Hilbert's projective metric \cite{lemmens2014birkhoff}. 

It is tempting to apply further change of variables $t\mapsto s:=T-t$, $\varphi(\bm{\xi},t) \mapsto p(\bm{\xi},s)$ proposed in \cite[Theorem 3]{caluya2021wasserstein} to \eqref{Coupled_PDE_first_order}, for transforming \eqref{Forward_PDE_First_Order}-\eqref{Backward_PDE_First_Order} into forward-forward PDEs as in \cite[equation (33)]{caluya2021wasserstein}. When possible, this strategy allows using a single FPK initial value problem (IVP) solver to set up a provably contractive fixed point recursion for computing the pair $\left(\hat{\varphi}_0,\varphi_{T}\right)$. In our case, the aforesaid mappings transform \eqref{Backward_PDE_First_Order} to
\begin{align}
&\!\!\!\!\!\!\frac{\partial p}{\partial s}=\nabla_{\bm{\xi}} \cdot\left(\!p \nabla_{\bm{\xi}} \tilde{V}\!\right)\!+\!\Delta_{\bm{\xi}} p\nonumber\\
&\!\!+\!\underbrace{p\left\langle\!\nabla_{\bm{\xi}} \tilde{V},(\bm{I}_{n}-\bm{\Upsilon}) \nabla_{\bm{\xi}} \tilde{V}\!\right\rangle\!+\!\left\langle\!\nabla_{\bm{\xi}} p,(\bm{I}_{n}-\bm{\Upsilon}) \nabla_{\bm{\xi}} \tilde{V}\!\right\rangle}_{\text {extra terms compared to  \cite[equation (33b)]{caluya2021wasserstein}}}, \label{backwardPDEwithEXtra_terms}
\end{align}
which has additional terms compared to  \cite[equation (33b)]{caluya2021wasserstein}. An interesting observation follows: \eqref{backwardPDEwithEXtra_terms} becomes the same forward FPK operator as in \eqref{Forward_PDE_First_Order} only if $\bm{\Upsilon}$ equals identity. Consequently, the Algorithm \textsc{ComputeFactorsSBP} proposed in  \cite[Sec. \RNum{5}.D]{{caluya2021wasserstein}} that uses a single FPK IVP solver, cannot be applied to our case. We need two different solvers for \eqref{Forward_PDE_First_Order} and \eqref{Backward_PDE_First_Order}.

To solve \eqref{Forward_PDE_First_Order}, we implement a modified form of the \textsc{ProxRecur} algorithm given in \cite[Sec. \RNum{3}.B]{caluya2019gradient} with the following distance functional, which is a weighted version of the squared 2-Wasserstein distance between a pair of joint PDFs $\tilde{\varrho},\tilde{\varrho}_{k-1}$, given by
\begin{align}
\!W_{\bm{\Upsilon}}^{2}\left(\tilde{\varrho},\tilde{\varrho}_{k-1}\right) := \!\!\underset{\pi\in\Pi\left(\tilde{\varrho},\tilde{\varrho}_{k-1}\right)}{\inf}\!\displaystyle\int_{\left(\prod_{i=1}^{n}[0,2\pi/\sigma_{i})\right)^{2}}\nonumber\\
\qquad\big\langle\bm{\theta}-\bm{\bar{\theta}},\bm{\Upsilon}^{-1}(\bm{\theta}-\bm{\bar{\theta}})\big\rangle\:\differential\pi(\bm{\theta},\bm{\bar{\theta}}),
\label{DefWassContinuous}    
\end{align}
where $\Pi\left(\tilde{\varrho},\tilde{\varrho}_{k-1}\right)$ is the set of joint probability measures supported on $\left(\prod_{i=1}^{n}[0,2\pi/\sigma_{i})\right)^{2}$, having finite second moments, with given marginal PDFs $\tilde{\varrho},\tilde{\varrho}_{k-1}$.

To solve \eqref{Backward_PDE_First_Order}, we employ the Feynman-Kac formula \cite{del2004feynman} as detailed in Sec. \ref{sec:feynman_kac}.

\subsection{Second Order Case}
In the second order Kuramoto model \eqref{SecondOrderKuramotoSDE}, the anisotropy in process noise directly affects the last $n$ components. Motivated by our treatment in the first order case, we now consider the invertible linear map 
\begin{align}
\begin{pmatrix}
\bm{\theta}\\
\bm{\omega}
\end{pmatrix}
\mapsto \begin{pmatrix}
\bm{\xi}\\
\bm{\eta}
\end{pmatrix}
:=\left(\bm{I}_{2}\otimes (\bm{M} \bm{S}^{-1})\right) \begin{pmatrix}
\bm{\theta}\\
\bm{\omega}
\end{pmatrix}
\label{invertiblemapSecond_order}
\end{align}
which by It$\hat{\text{o}}$'s Lemma \cite[Ch.4.2]{oksendal2013stochastic}, results in the following SDE for the transformed state vector $(\bm{\xi},\bm{\eta})^{\top}$: 
        {\small{\begin{equation}
             \begin{pmatrix}
                      \differential \bm{\xi}\\
                      \differential \bm{\eta}
             \end{pmatrix}=\begin{pmatrix}
                      \bm{\eta}\\
                     \bm{u} -\bm{\widetilde{\Upsilon}} \nabla_{\bm{\xi}} U(\bm{\xi})-\nabla_{\bm{\eta}} F(\bm{\eta})
             \end{pmatrix} \differential t + \begin{pmatrix}
                      \bm{0}_{n\times n}\\
                      \bm{I}_n
             \end{pmatrix} \differential \bm{w}
             \label{SDE_xi_eta}
         \end{equation}}}  
where $\bm{\widetilde{\Upsilon}}:=\left(\prod_{i=1}^{n}\sigma_{i}^{2}m_{i}^{-2}\right)\bm{M}\bm{S}^{-2}$, and the potentials 
{\small{\begin{align*}
U(\bm{\xi}):=&\left( \sum\limits_{i<j}k_{ij}\left(1-\cos\left(\frac{\sigma_{i}}{m_{i}}\xi_{i}-\frac{\sigma_{j}}{m_{j}}\xi_{j}-\varphi_{ij}\right)\right)-\right.\\
&\left.~~~~~\sum\limits_{i=1}^{n} \frac{\sigma_{i}}{m_{i}}P_{i}\xi_{i} \right)\left(\prod\limits_{i=1}^{n} \left( \frac{m_{i}}{\sigma_{i}}\right)^{\!2} \right),\\
F(\bm{\eta}):=&\frac{1}{2}\langle\bm{\eta},\bm{S}^{-1}\bm{\Gamma}\bm{\eta} \rangle.
\end{align*}}}
In this new state coordinate, the problem \eqref{FPKobjective} subject to \eqref{FPPDESecondOrder} and the endpoint PDF constraints, takes the form
\begin{subequations}
\begin{align}
&\underset{(\tilde{\rho}^{\bm{u}},\bm{u})}{\inf}\:\int_{0}^{T}\int_{\mathcal{X}}\|\bm{u}(\bm{\xi},\bm{\eta},t)\|_{2}^{2}\:\tilde{\rho}^{\bm{u}}(\bm{\xi},\bm{\eta},t)\:\differential \bm{\xi}\:\differential\bm{\eta}\:\differential t\\
&\begin{aligned}
	\frac{\partial \tilde{\rho}^{\bm{u}}}{\partial t} =&\nabla_{\bm{\eta}} \cdot\left(\tilde{\rho}^{\bm{u}} \left(-\bm{u} +\bm{\widetilde{\Upsilon}} \nabla_{\bm{\xi}} U(\bm{\xi})+\nabla_{\bm{\eta}} F(\bm{\eta}\right)\right)\\
	&-\langle \bm{\eta},\nabla_{\bm{\xi}} \tilde{\rho}^{\bm{u}} \rangle+\Delta_{\bm{\eta}} \tilde{\rho}^{\bm{u}},
	\end{aligned}\\
&{\small{\begin{aligned}
\tilde{\rho}^{\bm{u}}(\bm{\xi},\bm{\eta},0)&=\rho_{0}\!\left(\!\left(\!\bm{I}_{2}\otimes\bm{S}\bm{M}^{-1}\!\right)\!\begin{pmatrix}\bm{\xi}\\
\bm{\eta}
\end{pmatrix}\!\right)\!\!\left(\prod_{i=1}^{n} \frac{\sigma_i^{2}}{m_{i}^{2}}\right),\\
\tilde{\rho}^{\bm{u}}(\bm{\xi},\bm{\eta},T)&=\rho_{T}\!\left(\!\left(\!\bm{I}_{2}\otimes\bm{S}\bm{M}^{-1}\!\right)\!\begin{pmatrix}\bm{\xi}\\
\bm{\eta}
\end{pmatrix}\!\right)\!\!\left(\prod_{i=1}^{n} \frac{\sigma_i^{2}}{m_{i}^{2}}\right).
\end{aligned}}}
\end{align}
\label{opt_second_order}
\end{subequations}

Applying Proposition 1 and Theorem 2 of \cite{caluya2021wasserstein} to \eqref{opt_second_order}, we next derive a boundary-coupled system of linear PDEs akin to \eqref{Coupled_PDE_first_order}, for the Schr\"{o}dinger factors $(\varphi,\hat{\varphi})$, given by
{\small{\begin{subequations}
\begin{align}
&\begin{aligned}
\frac{\partial \hat{\varphi}}{\partial t}=-\langle \bm{\eta}, \nabla_{\bm{\xi}} \hat{\varphi} \rangle+ \nabla_{\bm{\eta}}\!\cdot\!\left(\! \hat{\varphi}(\bm{\widetilde{\Upsilon}} \nabla_{\bm{\xi}} \bm{U}(\bm{\xi})+\nabla_{\bm{\eta}} \bm{F}(\bm{\eta})) \!\right)\!+\!\Delta_{\bm{\eta}}\hat{\varphi} \end{aligned},\label{Forward_PDE_Second_Order}\\
&\begin{aligned}
\frac{\partial \varphi}{\partial t}=-\langle \bm{\eta}, \nabla_{\bm{\xi}} \varphi \rangle\!+\!\langle \bm{\widetilde{\Upsilon}} \nabla_{\bm{\xi}} \bm{U}(\bm{\xi})\!+\!\nabla_{\bm{\eta}} \bm{F}(\bm{\eta}), \nabla_{\bm{\eta}} \varphi \rangle \!-\! \Delta_{\bm{\eta}} \varphi\label{Backward_PDE_Second_Order}\end{aligned},\\
&\hat{\varphi}_{0}(\bm{\xi},\bm{\eta})\varphi_{0}(\bm{\xi},\bm{\eta}) \!=\! \rho_{0}\!\left(\!\!\left(\!\bm{I}_{2}\otimes\bm{S}\bm{M}^{-1}\!\right)\!\begin{pmatrix}\bm{\xi}\\
\bm{\eta}
\end{pmatrix}\!\right)\!\!\left(\prod_{i=1}^{n} \frac{\sigma_i^{2}}{m_{i}^{2}}\right),\label{CoupledBoundaryCondSecondOrder1}\\
&\hat{\varphi}_{T}(\bm{\xi},\bm{\eta})\varphi_{T}(\bm{\xi},\bm{\eta}) \!=\! \rho_{T}\!\left(\!\!\left(\!\bm{I}_{2}\otimes\bm{S}\bm{M}^{-1}\!\right)\!\begin{pmatrix}\bm{\xi}\\
\bm{\eta}
\end{pmatrix}\!\right)\!\!\left(\prod_{i=1}^{n} \frac{\sigma_i^{2}}{m_{i}^{2}}\right). \label{CoupledBoundaryCondSecondOrder2}
\end{align}
\label{Coupled_PDE_Second_order}
\end{subequations}}}

The optimal decision variables $(\tilde{\rho}^{\text{opt}},\bm{u}^{\text{opt}})$ for problem \eqref{opt_second_order} are obtained from the solution of \eqref{Coupled_PDE_Second_order} as
\begin{align}
 \tilde{\rho}^{\text{opt}}(\bm{\xi},\bm{\eta},t)&= \hat{\varphi}(\bm{\xi},\bm{\eta},t)\varphi(\bm{\xi},\bm{\eta},t),\nonumber\\
 \bm{u}^{\text{opt}}(\bm{\xi},\bm{\eta},t)&=\nabla_{\!\!{\tiny{{\begin{pmatrix}\bm{\xi}\\
 \bm{\eta}\end{pmatrix}}}}
 } \!\!\log \varphi (\bm{\xi},\bm{\eta},t).
 \label{opt_pair_2nd}
\end{align}
The optimally controlled joint state PDF $\rho^{\text{opt}}$ for \eqref{FPKobjective} in the second order case, is then obtained as $$\rho^{\text{opt}}(\bm{\theta},\bm{\omega},t)= \tilde{\rho}^{\text{opt}}\!\left(\!\left(\!\bm{I}_{2}\otimes\bm{M}\bm{S}^{-1}\!\right)\!\!\begin{pmatrix}\!\bm{\theta}\!\\
\!\bm{\omega}\!\end{pmatrix}\!,t\!\right)\!\left(\prod_{i=1}^{n}\frac{m_i^2}{\sigma_{i}^2}\right).$$ 
The optimal control in the original coordinates is $\left(\!\bm{I}_{2}\otimes\bm{S}\bm{M}^{-1}\!\right)\nabla_{\!\!{\tiny{{\begin{pmatrix}\bm{\theta}\\
 \bm{\omega}\end{pmatrix}}}}
 } \!\!\log\varphi\!\left(\!\left(\!\bm{I}_{2}\otimes\bm{M}\bm{S}^{-1}\!\right)\!\begin{pmatrix}\!\bm{\theta}\!\\
\!\bm{\omega}\!\end{pmatrix},t\!\right)$.

As in the first order case, our algorithmic approach (to be detailed in Sec. \ref{sec:overal_algorithm}) is to solve \eqref{Coupled_PDE_Second_order} via fixed point recursion over the pair $(\hat{\varphi}_{0},\varphi_{T})$ that is provably contractive w.r.t. the Hilbert's projective metric. In particular, to solve the backward Kolmogorov PDE \eqref{Backward_PDE_Second_Order}, we use the Feynman-Kac formula detailed in Sec. \ref{sec:feynman_kac}. The PDE \eqref{Forward_PDE_Second_Order} is the so-called \emph{kinetic Fokker-Planck equation} \cite[p. 40]{villani2009hypocoercivity}, and to solve the same, we propose a modified version of the proximal recursion proposed in \cite[Sec. \RNum{5}.B]{caluya2019gradient}. Our modification concerns with the distance functional in the proximal recursion, i.e., we consider the following analogue of \eqref{DefWassContinuous}:
\begin{align}
&\widetilde{W}_{h,\bm{\widetilde{\Upsilon}}}^{2}\left(\tilde{\varrho},\tilde{\varrho}_{k-1}\right) := \underset{\pi\in\Pi\left(\tilde{\varrho},\tilde{\varrho}_{k-1}\right)}{\inf}\!\displaystyle\int_{\left(\prod_{i=1}^{n}[0,2\pi m_{i}/\sigma_{i})\right)^{2} \times \mathbb{R}^{2n}}\nonumber\\
& \qquad\qquad\qquad\qquad s_{h,\bm{\widetilde{\Upsilon}}}\left(\bm{\xi},\bm{\eta},\bar{\bm{\xi}},\bar{\bm{\eta}}\right)\:\differential\pi\left(\bm{\xi},\bm{\eta},\bar{\bm{\xi}},\bar{\bm{\eta}}\right),
\label{DefWass_second_order}    
\end{align}
where $h>0$ is the step-size in proximal recursion, $\Pi\left(\tilde{\varrho},\tilde{\varrho}_{k-1}\right)$ is the set of joint probability measures over the product space $\left(\prod_{i=1}^{n}[0,2\pi m_{i}/\sigma_{i})\right)^{2} \times \mathbb{R}^{2n}$ that have finite second moments and marginal PDFs $\tilde{\varrho},\tilde{\varrho}_{k-1}$. The ``ground cost" in \eqref{DefWass_second_order} is
{\small{\begin{align}
     &s_{h,\bm{\widetilde{\Upsilon}}}\left(\bm{\xi},\bm{\eta},\bar{\bm{\xi}},\bar{\bm{\eta}}\right):=\nonumber\\
    &\bigg\langle\!\!\left(\bar{\bm{\eta}}-\bm{\eta}+h\bm{\widetilde{\Upsilon}} \nabla \bm{U}(\bm{\xi})\right),\bm{\widetilde{\Upsilon}}^{-1}\left(\bar{\bm{\eta}}-\bm{\eta}+h\bm{\widetilde{\Upsilon}} \nabla \bm{U}(\bm{\xi})\right)\!\!\bigg\rangle\nonumber\\
    & + 12 \bigg\langle\!\left(\frac{\bar{\bm{\xi}}-\bm{\xi}}{h}-\frac{\bar{\bm{\eta}}-\bm{\eta}}{h}\right),\bm{\widetilde{\Upsilon}}^{-1}\left(\frac{\bar{\bm{\xi}}-\bm{\xi}}{h}-\frac{\bar{\bm{\eta}}-\bm{\eta}}{h}\right)\!\bigg\rangle.
    \end{align}}}
In the next Section, we bring these ideas together to detail the algorithms for computing the optimal solutions in both the first and second order cases.

\section{Algorithms}\label{sec:algo}
In Sec. \ref{sec:prox_algo}, we first outline the proximal algorithm for solving the forward Kolmogorov PDEs \eqref{Forward_PDE_First_Order} and \eqref{Forward_PDE_Second_Order}. Then Sec. \ref{sec:feynman_kac} presents the Feynman–Kac algorithm for solving the backward Kolmogorov PDEs \eqref{Backward_PDE_First_Order} and \eqref{Backward_PDE_Second_Order}. Sec. \ref{sec:overal_algorithm} summarizes the  overall algorithm for solving \eqref{Coupled_PDE_first_order} and \eqref{Coupled_PDE_Second_order}.

\subsection{Proximal Algorithm}\label{sec:prox_algo}
For solving IVPs involving the forward Kolmogorov PDEs \eqref{Forward_PDE_First_Order} and \eqref{Forward_PDE_Second_Order}, we employ proximal recursions over the space of measurable positive functions over discrete time $t_{k-1}:=(k-1)h$ where the index $k\in\mathbb{N}$, and $h>0$ is (here constant) time step-size. These recursions are of the form
\begin{align}
\!\!\hat{\phi}_{k} \!=\! \prox_{h\Psi}^{d}\!\left(\!\hat{\phi}_{k-1}\!\right)\! := \!\arginf_{\hat{\phi}} \!\frac{1}{2}\!\left(\!d\!\left(\!\hat{\phi},\hat{\phi}_{k-1}\!\right)\!\right)^{\!\!2} \!\!+\! h\Psi\!\left(\!\hat{\phi}\!\right)  
\label{proxGeneric}    
\end{align}
where $\hat{\phi}_{k-1}(\cdot):=\hat{\phi}\left(\cdot,t_{k-1}\right)$, $d(\cdot,\cdot)$ is a distance-like functional, $\Psi$ is an energy-like functional, and $\hat{\phi}_{0}$ is suitable initial condition. The recursion \eqref{proxGeneric} reads as ``the proximal operator of the functional $h\Psi$ w.r.t. the distance $d$". The pair $(d,\Psi)$ is constructed in a way that the sequence of functions $\{\hat{\phi}_{k-1}\}_{k\in\mathbb{N}}$ generated by \eqref{proxGeneric} satisfies $\hat{\phi}_{k-1}(\cdot) \rightarrow \hat{\varphi}(\cdot,t)$ in $L^{1}\left(\mathcal{X}\right)$ as $h\downarrow 0$.

For \eqref{Forward_PDE_First_Order}, we set $d\equiv W_{\bm{\Upsilon}}$ given by \eqref{DefWassContinuous}, and $\Psi(\hat{\phi})\equiv \int_{\prod_{i=1}^{n}[0,2\pi/\sigma_i)} \left(\tilde{V} + \log\hat{\phi}\right)\hat{\phi}\:\differential\bm{\xi}$. For \eqref{Forward_PDE_Second_Order}, we set $d\equiv W_{h,\tilde{\bm{\Upsilon}}}$ given by \eqref{DefWass_second_order}, and $\Psi(\hat{\phi})\equiv \int_{\left(\prod_{i=1}^{n}[0,2\pi m_i/\sigma_i)\right)\times\mathbb{R}^{n}}\left(F + \log\hat{\phi}\right)\hat{\phi}\:\differential\bm{\xi}\:\differential\bm{\eta}$. For a discussion on the convergence guarantees and on implementation of these proximal updates via fixed point recursions, we refer the readers to \cite{caluya2019gradient}; see also \cite[Sec. V-B,C]{caluya2021wasserstein}. 

\subsection{Feynman-Kac Algorithm}\label{sec:feynman_kac}
For solving IVPs involving the backward Kolmogorov PDEs \eqref{Backward_PDE_First_Order} and \eqref{Backward_PDE_Second_Order}, we employ the Feynman-Kac path integral formulation \cite[Ch. 8.2]{oksendal2013stochastic}, \cite{yong1997relations}, \cite[Ch. 3.3]{yong1999stochastic}. The main idea is to solve the IVPs associated with \eqref{Backward_PDE_First_Order} and \eqref{Backward_PDE_Second_Order} using the sample path simulations of the corresponding uncontrolled \emph{forward} SDEs. We mention here that several works in stochastic control and learning \cite{pra1990markov,theodorou2010learning,pereira2020feynman,hawkins2021feynman} have leveraged the computational benefits of the Feynman-Kac approach. Specifically, the Feynman-Kac formula allows expressing the solution of backward PDE IVP
\begin{align*}
    &\frac{\partial \varphi}{\partial t}=\langle \nabla_{\tilde{\bm{x}}} \varphi , \bm{f}(\tilde{\bm{x}},t)\rangle +{\mathrm{trace}} \left(\bm{G}(\tilde{\bm{x}},t)\bm{G}(\tilde{\bm{x}},t)^{\top}\textbf{Hess}(\varphi)\right),\\
    &\varphi(\tilde{\bm{x}},t=T)=\varphi_{T}(\tilde{\bm{x}}),
\end{align*}
as the conditional expectation 
\begin{align}
\varphi(\tilde{\bm{x}},t) = \mathbb{E}\left[\varphi_{T}\left(\bm{x}(T)\right)\mid \bm{x}(t)=\tilde{\bm{x}}\right]
\label{ConditionalExpectation}
\end{align} 
where $\bm{x}(t)$ follows the It\^{o} diffusion $\differential\bm{x}(t) = \bm{f}(\bm{x},t)\differential t + \bm{G}(\bm{x},t)\differential\bm{w}$.

We use Algorithm \ref{alg:Feynman–Kac} to compute the Schr\"{o}dinger factor $\varphi$ as the conditional expectation \eqref{ConditionalExpectation} estimated from the forward SDE sample path simulations via the Euler-Maruyama scheme. The respective $\bm{f},\bm{G}$ for these sample path simulations correspond to those in the uncontrolled PDEs \eqref{Backward_PDE_First_Order} and \eqref{Backward_PDE_Second_Order}.

\begin{algorithm}
\caption{Feynman-Kac Algorithm for solving the backward PDE IVP at $t=\tau$}\label{alg:Feynman–Kac}
\begin{algorithmic}[1]
\Procedure{FeynmanKac}{$\varphi_{T}(\tilde{\bm{x}}_{T}),  \tilde{\bm{x}}_{T},T,\tilde{\bm{x}}_{\tau}, \tau, \bm{f},\bm{G},\newline N_{r},\text{nSample}, \text{dim},h,\bm{\lambda}$}
\State $\tilde{\bm{x}}_{r}\leftarrow [\bm{0}_{\text{nSample}\times \text{dim}\times N_{r}}]$ \Comment{{\small{initialize}}}
\State $\varphi_{r}\leftarrow [\bm{0}_{\text{nSample}\times N_{r}}]$
\State numSteps $\leftarrow (T-\tau)/h$ 
\For{$i$ $\leftarrow$ $1$ \text{to} $N_{r}$}
\State $\tilde{\bm{x}}_{\text{temp}}\leftarrow[\bm{x}_\tau,\bm{0}_{\text{nSample}\times \text{dim} \times \text{numSteps}}]$
 \For{$k$ $\leftarrow$ $1$ \text{to numSteps}}
\State$\tilde{\bm{x}}_{\text{temp}}\!(:,:,k+1\!)\!\!\leftarrow\!\!\tilde{\bm{x}}_{\text{temp}}(:,:,k)\!\!+\!\!h\bm{f}(\tilde{\bm{x}}(:,:,k),k)\!+\!\bm{G}(\tilde{\bm{x}}_{k},k)(\bm{w}_{k+1}\!-\!\bm{w}_{k})$ \!\!\!\Comment{{\small{Euler-Maruyama update}}}  
\EndFor
\!\State $\tilde{\bm{x}}_{r}(:,:,i)\leftarrow \tilde{\bm{x}}_{\text{temp}}(:,:,\text{numSteps}+1)$
\!\State $\varphi_{r}(:,i)\!\!\leftarrow\!\text{ElasticNet} (\varphi_{T}(\tilde{\bm{x}}_T),\tilde{\bm{x}}_{T},\tilde{\bm{x}}_{r}(:,:,i),\bm{\lambda})$ \label{ElasticNetReg}
\EndFor
\State \Return  $\varphi(\tilde{\bm{x}}_\tau,\tau)\leftarrow\frac{1}{N_{r}} \sum\limits_{i=1}^{N_r} \varphi_{r}(:,i)$
\EndProcedure
\end{algorithmic}
\end{algorithm}

In Algorithm \ref{alg:Feynman–Kac}, $\bm{\lambda}:=(\lambda_1,\lambda_2)\in\mathbb{R}^{2}_{>0}$ is a regularizing parameter vector. In line \ref{ElasticNetReg} of Algorithm \ref{alg:Feynman–Kac}, we implement an elastic net regression \cite{zou2005regularization}, referred to as ``ElasticNet", with $\lambda_{1},\lambda_{2}$ being the regularizing weights for the $\ell_1$ and $\ell_2$ norms, respectively. ElasticNet approximates the value of $\varphi_{T}$ at $\tilde{\bm{x}}_{r}(:,:,i)$ from the known boundary values $\varphi_{T}(\tilde{\bm{x}}_T)$ and the propagated samples $\tilde{\bm{x}}_{T}$. For the simulation results reported in Sec. \ref{sec:numericalsimulation}, the ElasticNet computes a degree three polynomial approximant in the transformed state co-ordinates. We use the Alternating Direction Method of Multipliers (ADMM) algorithm \cite[Ch. 6]{boyd2011distributed} to implement the elastic net regression. We estimate \eqref{ConditionalExpectation} as an empirical average (line 13 of Algorithm \ref{alg:Feynman–Kac}) of the approximated $\varphi$ at time $t=T$ over $N_{r}$ sample paths. The parameters nSample and dim in Algorithm \ref{alg:Feynman–Kac} denote the number of samples and the state dimension ($n$ for first order Kuramoto, $2n$ for second order Kuramoto case), respectively.

\subsection{Overall Algorithm}\label{sec:overal_algorithm}
Bringing together the ideas from Sec. \ref{sec:prox_algo} and \ref{sec:feynman_kac}, we now outline the overall algorithm to solve \eqref{Coupled_PDE_first_order} or \eqref{Coupled_PDE_Second_order}. To keep notations succinct, let us use $\tilde{\bm{x}}$ as the appropriate transformed state, i.e., $\tilde{\bm{x}}\equiv\bm{\xi}$ for the first order case, and $\tilde{\bm{x}}\equiv(\bm{\xi},\bm{\eta})$ for the second order case. We perform a fixed point recursion over the pair $\left(\hat{\varphi}_0(\tilde{\bm{x}}),\varphi_{T}(\tilde{\bm{x}})\right)$ as follows.

\textbf{Step 1.} Initialize arbitrary $\varphi_0(\tilde{\bm{x}})$ everywhere positive.

\textbf{Step 2.} Compute $\hat{\varphi}_0(\tilde{\bm{x}})=\tilde{\rho}^{\bm{u}}(\tilde{\bm{x}},0)/ \varphi_0(\tilde{\bm{x}})$. 

\textbf{Step 3.} Using \eqref{proxGeneric}, solve IVP \eqref{Forward_PDE_First_Order} or \eqref{Forward_PDE_Second_Order}  till $t=T$ to obtain $\hat{\varphi}_T(\tilde{\bm{x}})$.

\textbf{Step 4.} Compute $ \varphi_T(\tilde{\bm{x}})=\tilde{\rho}^{\bm{u}}(\tilde{\bm{x}})/\hat{\varphi}_T(\tilde{\bm{x}})$. 

\textbf{Step 5.} Use Algorithm \ref{alg:Feynman–Kac} to calculate $\varphi_0(\tilde{\bm{x}}):=\varphi(\tilde{\bm{x}},0)$ for \eqref{Backward_PDE_First_Order} or  \eqref{Backward_PDE_Second_Order}.

\textbf{Step 6.} Repeat until the pair $\left(\hat{\varphi}_0(\tilde{\bm{x}}),\varphi_{T}(\tilde{\bm{x}})\right)$ has converged w.r.t. the Hilbert's projective metric \cite{lemmens2014birkhoff}.

\textbf{Step 7.} Compute the Schr\"{o}dinger factors $\left(\hat{\varphi}(\tilde{\bm{x}},t),\varphi(\tilde{\bm{x}},t)\right)$ using the IVPs \eqref{Coupled_PDE_first_order} and \eqref{Coupled_PDE_Second_order}.

\textbf{Step 8.} Use $\left(\hat{\varphi}(\tilde{\bm{x}},t),\varphi(\tilde{\bm{x}},t)\right)$ from \textbf{Step 7} to compute the pair $(\tilde{\rho}^{\text{opt}},\bm{u}^{\text{opt}})$ from \eqref{opt_pair} or \eqref{opt_pair_2nd}.

\textbf{Step 9.} Bring back the optimal joint state PDF and the optimal control to the original coordinates, i.e., to $\bm{\theta}$ for the first order, and to $(\bm{\theta},\bm{\omega})$ for the second order case.

Since the fixed point recursion over the function pair $\left(\hat{\varphi}_0,\varphi_{T}\right)$ is contractive \cite{chen2016entropic} in Hilbert's projective metric, the above nine step algorithm is \emph{guaranteed to converge to a unique solution}.


\section{Numerical Simulations}\label{sec:numericalsimulation}
\subsubsection{First order case}\label{sec:sim_1st_order}
We consider an instance of \eqref{FirstOrderKuramotoSDE} with $n=2$ oscillators, i.e., $\bm{\theta} \in \mathcal{X}=\mathbb{T}^2$. We generated the following parameters uniformly random from the respective intervals: $P_i \in [0,10]$, $\sigma_i \in [1,5]$ for $i=1,2$, and $k_{12} \in [0.7,1.2]$, $\varphi_{12}\in [0,\frac{\pi}{2})$.

\begin{figure}[htbp!]
    \centering
\begin{subfigure}{.28\linewidth}
        \includegraphics[scale=0.3]{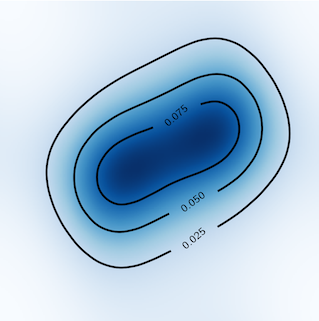}
        \caption{$t=0$}
    \end{subfigure}
    \hskip3em
    \begin{subfigure}{.28\linewidth}
        \includegraphics[scale=0.3]{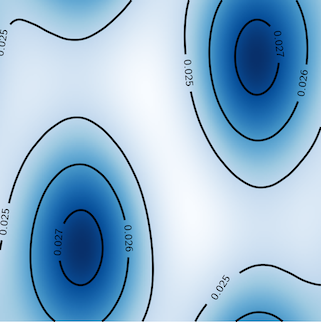}
        \caption{$t=1$}
    \end{subfigure}
\caption{%
{\small{Endpoint von Mises $\bm{\theta}$ PDFs over $\mathbb{T}^{2}$.}}}%
\label{fig:rho_0_T}
\vspace*{-0.15in}
\end{figure}

We set the final time $T=1$, and $\rho_{0},\rho_{T}$ as in \eqref{endpointPDFs} (see Fig. \ref{fig:rho_0_T}) with $\bm{\kappa}_{0}=\left(1,1\right)^{\!\top}$, $\bm{\kappa}_{T}=\left(0.01,0.01\right)^{\!\top}$, $\bm{m}_{0}=\left(\pi,\pi\right)^{\!\top}$, $\bm{m}_{T}=\left(0,0\right)^{\!\top}$, $\bm{\Lambda}_{0}=\begin{bmatrix} 0& 1\\1&0\end{bmatrix}$, $\bm{\Lambda}_{T}=0.1\bm{\Lambda}_{0}$. We solve \eqref{Coupled_PDE_first_order} following the steps in Section \ref{sec:overal_algorithm}. Specifically, we solve the backward PDE \eqref{Backward_PDE_First_Order} via Algorithm \ref{alg:Feynman–Kac} with parameters $N_r=100, h=0.1, \text{nSample}=441, \lambda_1=\lambda_2=0.01$. To solve the forward PDE \eqref{Forward_PDE_First_Order}, we used the PROXRECUR algorithm from \cite[Sec. III-B.1]{caluya2019gradient} with algorithmic parameters $\varepsilon=1,\beta=0.1,\delta=0.1,L=300$ together with the modifications mentioned in Sec. \ref{sec:prox_algo}.

\begin{figure*}[htbp!]
    \centering
      \begin{subfigure}{0.8\paperwidth}
        \includegraphics[width=0.8\paperwidth]{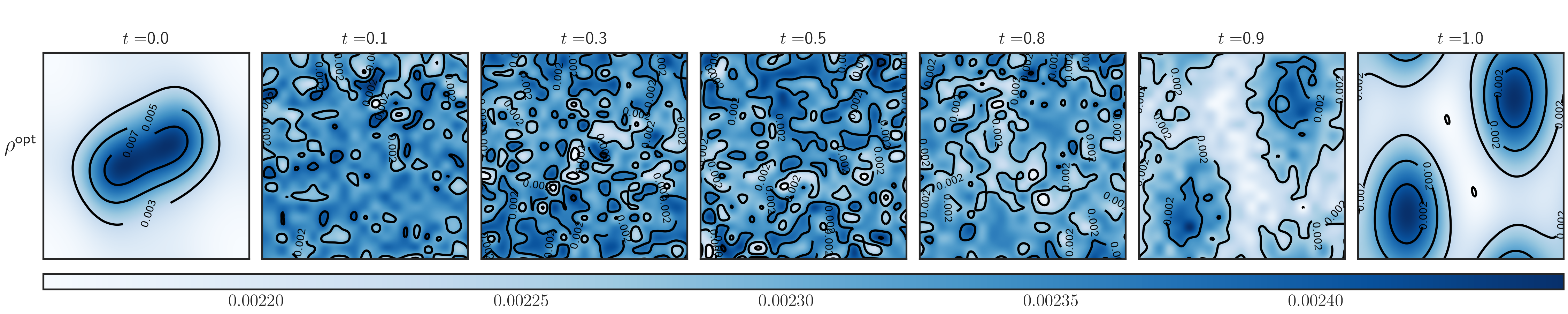}
          \caption{{\small{Contour plots of the optimally controlled state PDFs $\rho^{\text{opt}}(\bm{\theta},t)$ over $\mathbb{T}^{2}$. Each subplot is a snapshot in time $t\in [0,1]$.}}}
          \label{fig:opt_rho_1st_order}
      \end{subfigure}
        \hfill
      \begin{subfigure}{0.8\paperwidth}
        \includegraphics[width=0.8\paperwidth]{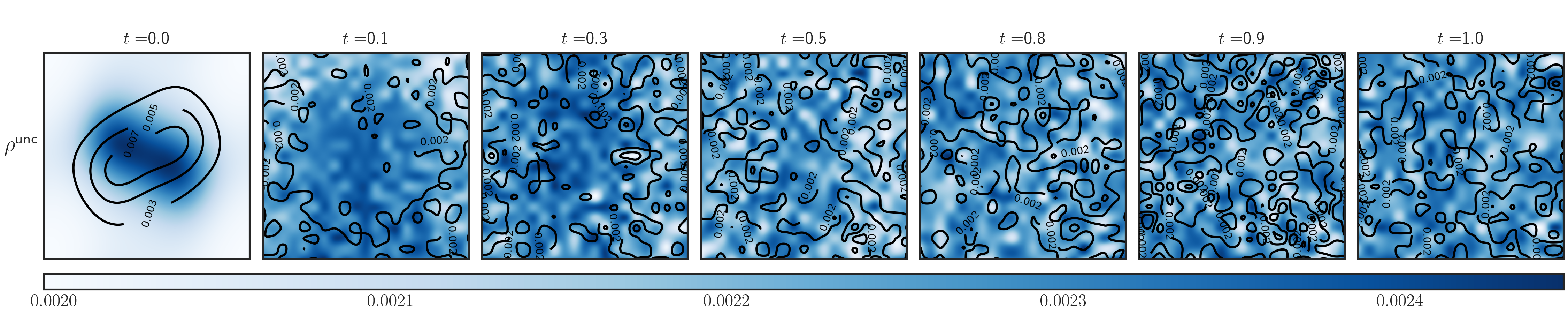}
          \caption{Contour plots of the uncontrolled state PDFs $\rho^{\text{unc}}(\bm{\theta},t)$ over $\mathbb{T}^{2}$. Each subplot is a snapshot in time $t\in [0,1]$.}
          \label{fig:unc_rho_1st_order}
      \end{subfigure}
      \hfill
      \begin{subfigure}{0.8\paperwidth}
        \includegraphics[width=0.8\paperwidth]{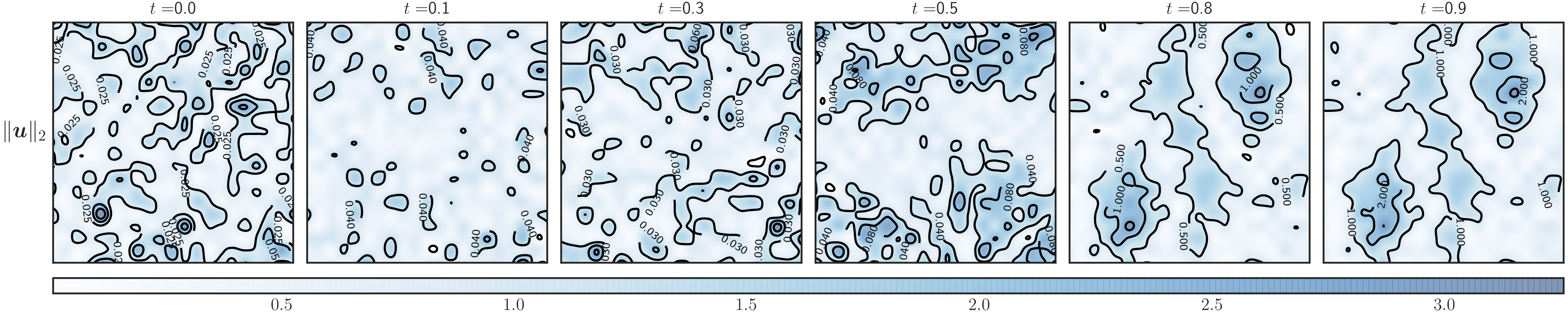}
          \caption{Contour plots of the 2-norm magnitude of the optimal control over $\mathbb{T}^{2}$. Each subplot is a snapshot in time $t\in [0,1]$.}
          \label{fig:opt_control_1st_order}
      \end{subfigure}
\caption{%
{\small{Simulation results for the optimal PDF steering for the \emph{first order} Kuramoto oscillators over $t\in[0,1]$. The color denotes the value of the plotted variable; see colorbar (dark hue = high, light hue = low).}}}
\label{fig:1st_order}
\vspace*{-0.2in}
\end{figure*}

Fig. \ref{fig:opt_rho_1st_order} shows the snapshots of the \emph{optimally controlled joint} $\rho^{\text{opt}}(\bm{\theta},t)$ steering $\rho_0$ to $\rho_{T}$ over time horizon $[0,1]$. Fig. \ref{fig:unc_rho_1st_order} shows the snapshots of the uncontrolled joint $\rho^{\text{unc}}(\bm{\theta},t)$ from the same $\rho_{0}$. The snapshots of the magnitude of optimal control are depicted in Fig. \ref{fig:opt_control_1st_order}.

\subsubsection{Second order case}\label{sec:sim_2nd_order}
We next consider an instance of \eqref{SecondOrderKuramotoSDE} with $n=2$ oscillators, i.e., $(\bm{\theta}, \bm{\omega})\in \mathcal{X} = \mathbb{T}^2 \times \mathbb{R}^2$. We set $T=1$, and use $\{P_i,\sigma_{i}\}_{i=1,2}$, $k_{12}$, $\varphi_{12}$ as in the first order case above. We consider the initial joint PDF $\rho_{0}(\bm{\theta},\bm{\omega})\equiv\overline{\rho}_{0}(\bm{\theta})\times\text{Unif}\left([0,0.2]^{2} \right)$, and the terminal joint PDF  $\rho_{T}(\bm{\theta},\bm{\omega})\equiv\overline{\rho}_{T}(\bm{\theta})\times \text{Unif}\left([0,0.2]^{2} \right)$ where the $\bm{\theta}$ marginals $\overline{\rho}_{0},\overline{\rho}_{T}$ are identical to $\rho_{0},\rho_{T}$ in the first order case, and Unif$(\cdot)$ denotes the uniform PDF. In other words, the endpoint joint PDFs $\rho_{0},\rho_{T}$ are supported on the compact set $\mathbb{T}^{2}\times[0,0.2]^2$. 
\begin{figure*}
      \begin{subfigure}{0.8\paperwidth}
        \includegraphics[width=0.8\paperwidth]{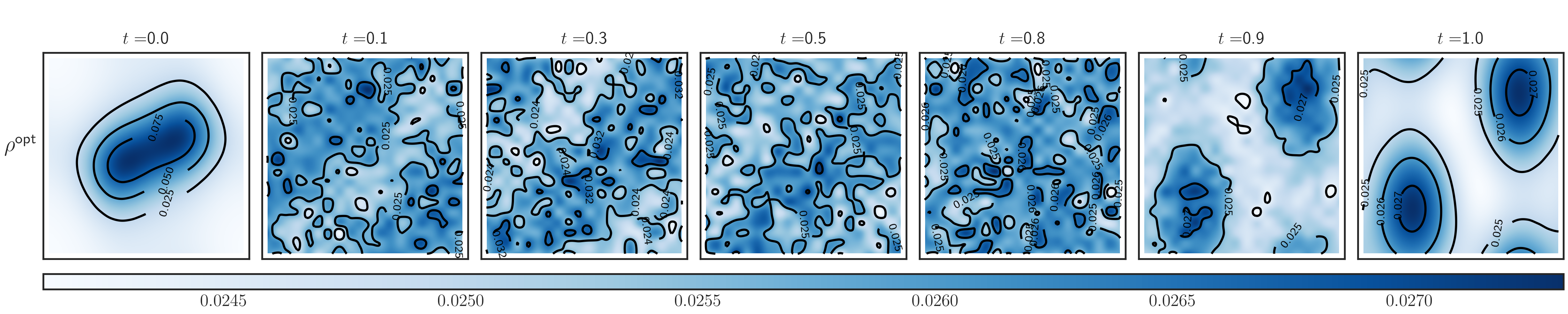}
          \caption{Contour plots of $\bm{\theta}$ marginals of the optimally controlled joints $\rho^{\text{opt}}(\bm{\theta},\bm{\omega},t)$ over $\mathbb{T}^{2}$. Each subplot is a snapshot in time $t\in [0,1]$.}
          \label{fig:opt_rho_2nd_order}
      \end{subfigure}
        \hfill
        \newline
      \begin{subfigure}{0.8\paperwidth}
        \includegraphics[width=0.8\paperwidth]{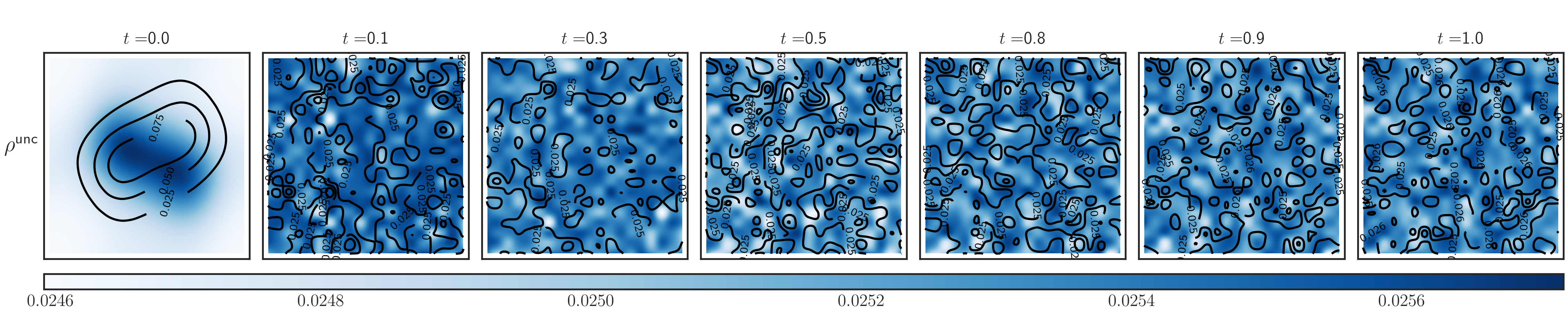}
          \caption{Contour plots of $\bm{\theta}$ marginals of the uncontrolled joints $\rho^{\text{unc}}(\bm{\theta},\bm{\omega},t)$ over $\mathbb{T}^{2}$. Each subplot is a snapshot in time $t\in [0,1]$.}
          \label{fig:unc_rho_2nd_order}
      \end{subfigure}
      \hfill
      \newline
      \begin{subfigure}{0.8\paperwidth}
        \includegraphics[width=0.8\paperwidth]{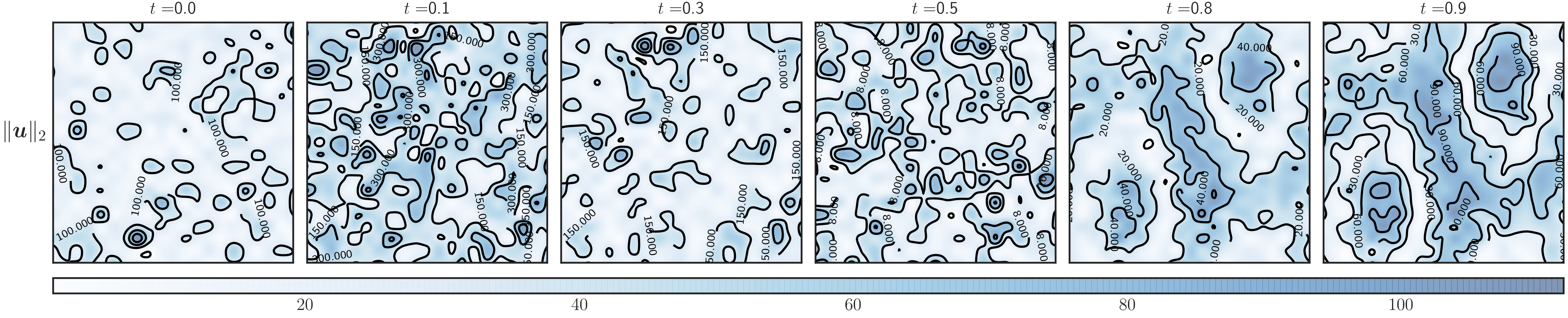}
          \caption{Contour plots of the 2-norm magnitude of the optimal control over $\mathbb{T}^{2}$. Each subplot is a snapshot in time $t\in [0,1]$.}
          \label{fig:opt_control_2nd_order}
      \end{subfigure}
\caption{
{\small{Simulation results for the optimal PDF steering for the \emph{second order} Kuramoto oscillators over $t\in[0,1]$. The color denotes the value of the plotted variable; see colorbar (dark hue = high, light hue = low).}}}
\label{fig:2nd_order}
\vspace*{-0.2in}
\end{figure*}

We solve \eqref{Coupled_PDE_Second_order} using the same computational set up as in the subsection above except that the PROXRECUR algorithm \cite[Sec. III-B.1]{caluya2019gradient} for solving the forward PDE \eqref{Forward_PDE_Second_Order} is suitably modified as mentioned in Sec. \ref{sec:prox_algo}.   

\begin{figure}[htbp!]
    \centering    
    \includegraphics[width=6.3cm]{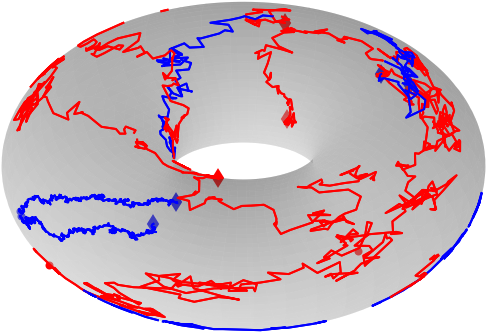}
\caption{%
{\small{The optimally controlled first order (\emph{in red}) and second order (\emph{in blue}) Kuramoto sample paths on $\mathbb{T}^{2}$ for the numerical simulation in Sec. \ref{sec:numericalsimulation}. The \emph{circled} and \emph{diamond markers} denote the initial and terminal angular coordinates, respectively.}}}%
\label{fig.optimallycontrolledsamplepaths}
\vspace*{-0.15in}
\end{figure}

Fig. \ref{fig:opt_rho_2nd_order} shows the snapshots of the $\bm{\theta}$ marginals of the optimally controlled joints $\rho^{\text{opt}}(\bm{\theta},\bm{\omega},t)$. Fig. \ref{fig:unc_rho_2nd_order} shows the $\bm{\theta}$ marginal snapshots of the uncontrolled joints. The snapshots of the magnitude of optimal control are depicted in Fig. \ref{fig:opt_control_2nd_order}. A comparison of Figs. \ref{fig:opt_control_1st_order} and \ref{fig:opt_control_2nd_order} reveals that in the second order case, the prior dynamics being mixed conservative-dissipative, the optimal control entails forcing that is about two orders of magnitude above the same for the first order case. Fig. \ref{fig.optimallycontrolledsamplepaths} shows four optimally controlled sample paths on $\mathbb{T}^{2}$ for the first order case (\emph{in red}) and another four for the second order case (\emph{in blue}).

\subsubsection{Order parameter}
\begin{figure*}[tpb!]
\centering
\parbox{8cm}{
\includegraphics[width=8cm]{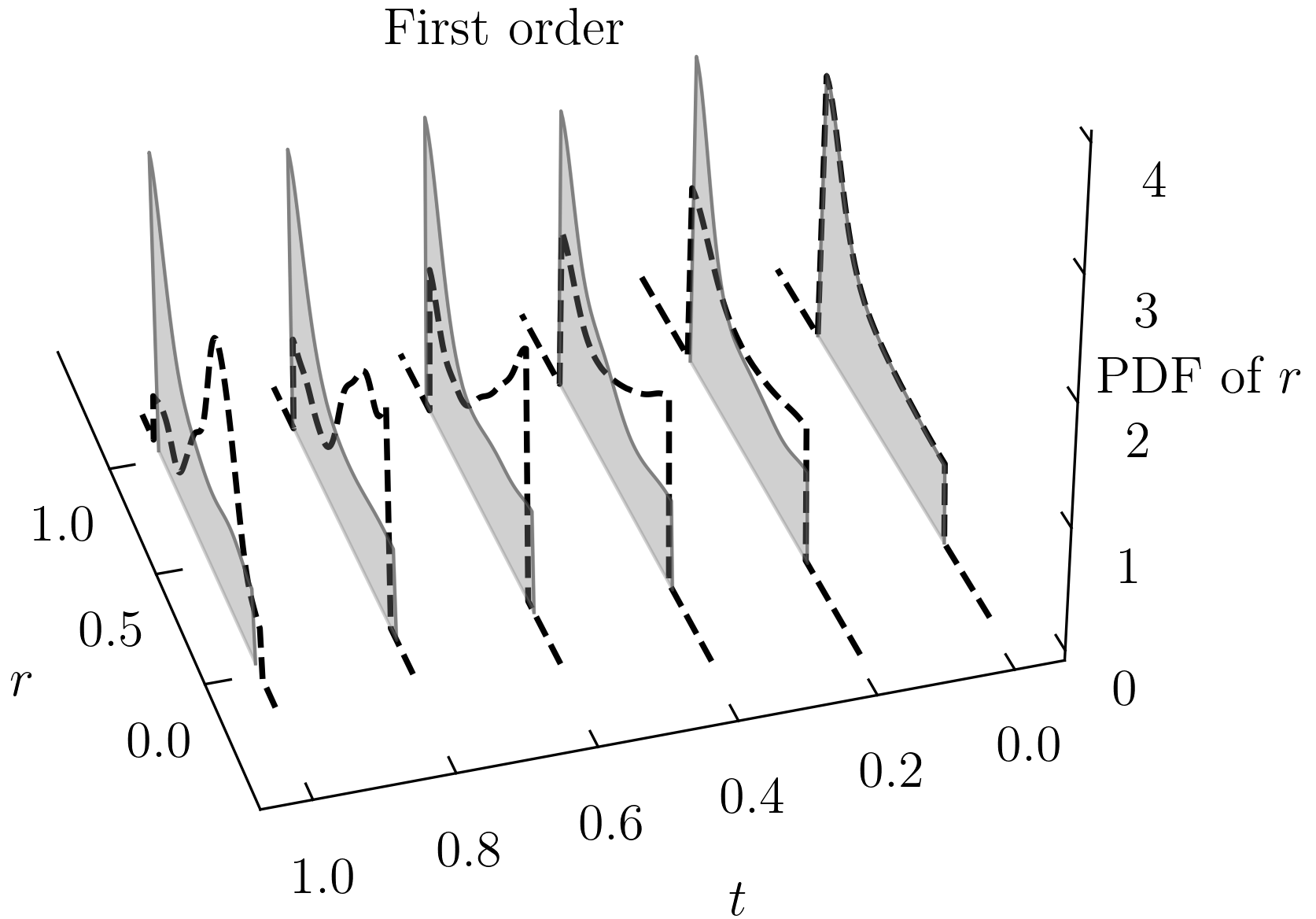}
}
\quad
\begin{minipage}{8cm}
\includegraphics[width=8cm]{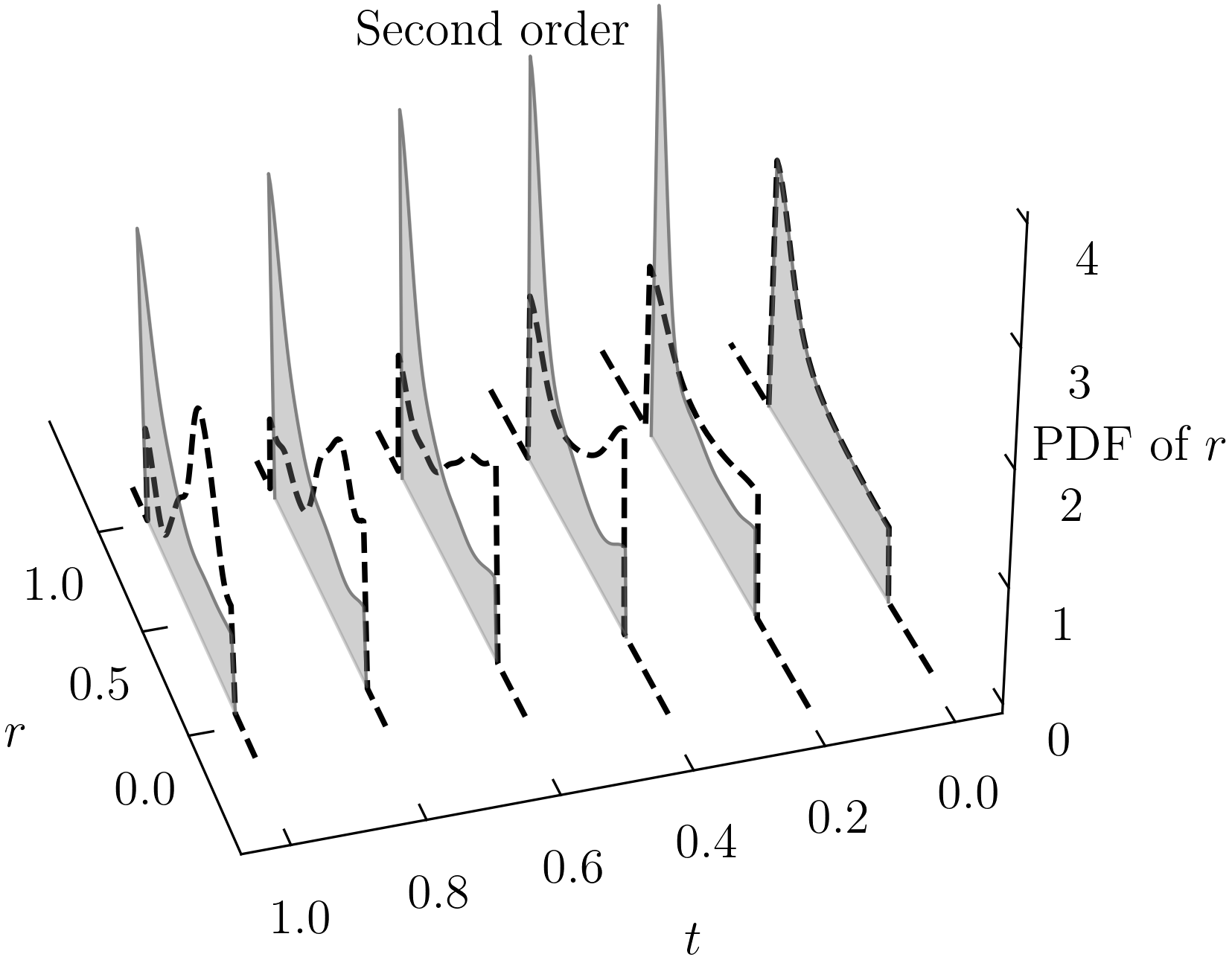}
\end{minipage}
\caption{%
{\small{PDFs of $r$ for the numerical simulation in Sec. \ref{sec:numericalsimulation}. The solid curves with grey filled areas are for the \emph{optimally controlled} $r$ PDFs. The dashed curves are for the \emph{uncontrolled} $r$ PDFs.}}}%
\label{fig.order parameters}
\vspace*{-0.15in}	
\end{figure*}

In the coupled oscillator context, a measure of synchronization, or lack thereof, is the order parameter $r:=\frac{1}{n}\vert\sum_{j=1}^{n} \exp(\iota \theta_{j})\vert \in [0,1]$ where $\iota:=\sqrt{-1}$; see e.g., \cite[Sec. 3.2]{strogatz2000kuramoto}.
For instance, $r=0$ implies lack of synchrony, and $r=1$ implies synchronized motion in the state space. Fig. \ref{fig.order parameters} shows the snapshots of the order parameter PDFs (solid lines with grey filled areas) under optimal control for the aforesaid numerical simulation and the order parameter PDFs for the uncontrolled cases (dashed lines). 

As the optimal control steers the stochastic state $\bm{\theta}$ from unimodal to bimodal, the optimally controlled $r$ PDFs (solid lines with grey filled areas) in Fig. \ref{fig.order parameters} slightly flatten over this transfer horizon and develop a secondary peak around $r=0.5$. The \emph{uncontrolled} $r$ PDFs (dashed curves) in Fig. \ref{fig.order parameters} show that as time progresses, the uncontrolled order parameter concentrates around $r=0$ indicating mixing/disorder with high probability, which is indeed consistent with the contour plots in Figs. \ref{fig:unc_rho_1st_order}-\ref{fig:unc_rho_2nd_order}.

\section{Conclusions}\label{sec:conclusions}
This paper proposes an algorithmic framework to solve the problem of minimum effort steering of the joint state PDF for a finite population of coupled noisy nonuniform Kuramoto oscillators subject to hard deadline and endpoint PDF constraints. This is an atypical stochastic control problem that is relevant to engineering applications such as active uncertainty management in a network of synchronous machines and loads. The problem differs from the infinite population a.k.a. the mean-field limit more commonly studied in the physics literature. The present work points out technical difficulties in applying some change-of-variables proposed in related literature for a class of dynamical nonlinearities, and alleviates the same by combining certain measure-valued proximal recursions with the Feynman-Kac path integral computation. Illustrative numerical simulations are given to highlight the results. The formulations and results presented herein, should be of broad interest to the researchers in systems, control and mathematical physics.


\bibliographystyle{IEEEtran}
\bibliography{References.bib}

\appendix
\subsection{ADMM for the  Elastic Net Regression}

In line \ref{ElasticNetReg} of Algorithm \ref{alg:Feynman–Kac}, we approximate the value of $\varphi_{T}$ at $\tilde{\bm{x}}_{r}(:,:,i)$ from the known boundary values $\varphi_{T}(\tilde{\bm{x}}_T)$ by computing a degree 3 polynomial approximation in 2 vector variables ($\bm{\xi}_1$, $\bm{\xi}_2$) for the first order case, and in 4 vector variables ($\bm{\xi}_1$,$\bm{\eta}_1$,$\bm{\xi}_2$,$\bm{\eta}_2$) for the second order case. 

The elastic net objective comprises of a squared 2 norm error for data fidelity, and additional 1 and squared 2 norm regularizations on the decision variable for parsimony. The decision variable in this case is the monomial coefficient vector $\bm{c}$ for the approximating multivariate cubic polynomial.

Specifically, given $\bm{b}:=\varphi_{T}(\tilde{\bm{x}}_T)$, the elastic net solves
 \begin{align}
     \min_{\bm{c}} \frac{1}{2}\|\bm{A} \bm{c}- \bm{b} \|_{2}^{2}+\lambda_{1}\|\bm{c}\|_{1}+\lambda_{2}\|\bm{c}\|_{2}^{2},
     \label{equ_elasnet}
 \end{align}
 where $\bm{A}$ is the data-dependent regression matrix, and $\lambda_1,\lambda_2>0$ are the regularizing coefficients. The hyperparameters $(\lambda_1,\lambda_2)\in\mathbb{R}^2_{>0}$ are typically optimized via implicit differentiation \cite{bertrand2020implicit} or Bayesian optimization \cite{snoek2012practical}.

 Letting  $\tilde{\bm{A}}:=\begin{bmatrix}\bm{A}\\\sqrt{\lambda_2}\bm{I}\end{bmatrix}$ and $\bm{\tilde{b}}:=\begin{bmatrix}\bm{b}\\ \bm{0}\end{bmatrix}$, we re-write \eqref{equ_elasnet} in the lasso \cite{tibshirani1996regression} form:
 \begin{align}
   \min_{\bm{c}}  \frac{1}{2}\|\tilde{\bm{A}} \bm{c}- \bm{\tilde{b}} \|_{2}^{2}+\lambda_{1}\|\bm{c}\|_{1},
 \end{align}
 and apply the corresponding (unscaled) ADMM recursion (see e.g., \cite[Ch. 6.4]{boyd2011distributed}) 
 \begin{align*}
 &\bm{c}_{j+1}=\left(\tilde{\bm{A}}^{\top} \tilde{\bm{A}}+r I\right)^{-1}\left(\tilde{\bm{A}}^{\top} \tilde{\bm{b}} + r \bm{z}_{j}-\bm{\nu}_{j}\right), \\
 &\bm{z}_{j+1}=S_{\lambda_1 / r}\left(\bm{c}_{j+1}+\bm{\nu}_{j} / r\right), \\
 &\bm{\nu}_{j+1}=\bm{\nu}_{j} + r\left(\bm{c}_{j+1}-\bm{z}_{j+1}\right),
 \end{align*}
where $j\in\mathbb{N}_0$ denotes the ADMM iteration index, the augmented Lagrangian regularizer $r=0.005$, and the soft thresholding operator $S_a(\cdot)$ for $a>0$ is defined elementwise for a vector argument as
\begin{align}
    \left(S_{a}(\bm{w})\right)_i := \begin{cases}w_{i}-a & \text { if } \quad w_{i} \geqslant a, \\ 0 & \text { if }-a \leqslant w_{i} \leqslant+a, \\ w_{i}+a & \text { if } \quad w_{i} \leqslant-a.\end{cases}
\end{align}

\end{document}